\documentstyle{amsart}

\def\E{{\hbox{\bf E}}}

 at 10 true pt

\def\be#1{ \begin{equation}\label{#1} }
\def\bas{\begin{align*}}
\def\eas{\end{align*}}
\def\bi{\begin{itemize}}
\def\ei{\end{itemize}}

\newenvironment{proof}{\noindent {\bf Proof} }{\endprf\par}
\def \endprf{\hfill  {\vrule height6pt width6pt depth0pt}\medskip}
\def\emph#1{{\it #1}}
\def\textbf#1{{\bf #1}}

\def\BZ{{\mathbf Z}}

\def\ep{{\epsilon}}
\theoremstyle{plain}
  \newtheorem{theorem}[subsection]{Theorem}

  \newtheorem{fact}[subsection]{Fact}
  \newtheorem{lemma}[subsection]{Lemma}
  \newtheorem{corollary}[subsection]{Corollary}
\theoremstyle{remark}
  \newtheorem{remark}[subsection]{Remark}

\theoremstyle{definition}
  \newtheorem{definition}[subsection]{Definition}
\include{psfig}
\begin{document}
\title[Structure of large incomplete sets in abelian groups ]
{ Structure of large incomplete sets in abelian groups}
%\author{Terence Tao}
%\address{Department of Mathematics, UCLA, Los Angeles CA 90095-1555}
%\email{tao@@math.ucla.edu}
%\thanks{T. Tao is a Clay Prize Fellow and is supported by a grant from the Packard Foundation.}

\author{Van H. Vu}
\address{Department of Mathematics, Rutgers, Piscataway, NJ 08854-8019}
\email{vanvu@@math.rutgers.edu}
\thanks{V. Vu is an A. Sloan Fellow and is supported by an NSF Career Grant.}
\begin{abstract}
Let $G$ be a finite abelian group and $A$ be a  subset of $G$. We
say that $A$ is complete if every element of $G$ can be
represented as a sum of different elements of $A$. In this paper,
we study the following question: \vskip2mm \centerline{\it What is
the structure of a large incomplete set ? } \vskip2mm The typical
answer is that such a set is essentially contained in a maximal
subgroup. As a by-product, we obtain a new proof for several
earlier results.
\end{abstract}

\maketitle
\section{Introduction}
Let $G$ be an abelian group and  $A$ be a subset of $G$. We use
$S_A$ to denote the  collection of all subset sums of $A$
$$S_A:= \{\sum_{a \in B} a| B \subset A, |B| < \infty \}. $$
\noindent We will keep this notation when $A$ is sequence of (not
necessarily different) elements of $A$. In this case $S_A$ is the
collection of all subsequence sums of $A$. $Z_n$ denotes the
cyclic group of order $n$. \vskip2mm
 {\it Example.} Take $G= Z_{11}$. If  $A$ is the subset
$\{1,2,3 \}$, then $S_A=\{1,2,3,4,5,6\}$. If $A$ is the sequence
$\{1,1,3\}$, then $S_A= \{1,2,3,4,5 \}$. \vskip2mm

 Following Erd\H os \cite{erd}, we say that $A$ is
{\it complete} if $S_A=G$ and {\it incomplete } otherwise. If $G$
is finite, the {\it critical} number of $G$, $c(G)$, is the
smallest integer $m$ such that any subset $A \subset G \backslash
\{0\}$ with size $m$ is complete. This parameter has been studied
for a long time and its exact value  is known for most groups.
\begin{theorem} \label{theorem:old} Let $G$ be a finite abelian
group of order $n=ph$, where $p$ is the smallest prime divisor of
$G$. Then the following holds
\begin{itemize}
\item If $p=2$ and $h \ge 5$ or $G= Z_2 \oplus Z_2 \oplus Z_2$,
then $c(G)= h$. If $p=2$ and $h \le 4$ and $G \neq Z_2 \oplus Z_2
\oplus Z_2$, then $c(G)= h+1$. \item If $h$ is a prime, then
$p+h-2 \le c(G) \le p+h-1$. Furthermore, if $h=p \ge 3$ or $h \ge
2p+1$, then $c(G)= p+h-2$. \item If $p \ge 3$ and $h$ is
composite, then $c(G) = p+h-2$.
\end{itemize}
\end{theorem}

The first statement is due to Diderrich and Mann \cite{DM}. The
second combines results of Mann and Wou \cite{MW} ( who studied
the case $h=p$) and Didderich \cite{D1}. The last statement has
been known as Didderich conjecture, posed in  \cite{D1} and was
proved by Gao and Hamidoune \cite{GH}, more than twenty years
later.

In this paper, we would like  to study the following question
\vskip2mm \centerline{\it What is the structure of a relatively
large incomplete set ? } \vskip2mm

Technically speaking, we would like to have a characterization for
incomplete sets of relatively large size. Such a characterization
has been obtained recently in \cite{GHLS} for sets of size at
least $n/(p+2)$. In this paper, we will be able to treat much
smaller sets. (In fact, our assumption on "relatively large" is
almost sharp; see Theorem \ref{theorem:dir2}.) The method used in
our proofs is different from those used in previous papers. As a
by-product, one obtains a new proof for a good portion of Theorem
\ref{theorem:old}, including a new proof for  Didderich's
conjecture for large $n$ (see the remarks following Theorem
\ref{theorem:dir1}).

{\it Notation.} $<A>$ denotes the subgroup generated by $A$.
$\E(X)$ denotes the expectation of a random variable $X$. All
logarithms have natural base, if not specified otherwise.

\section {The characterization of large incomplete sets}
Let us start by a simple fact, whose proof is left as an exercise.
\begin{fact} \label{fact1} Let $p$ be a prime and  $A$ be a sequence of $p-1$
non-zero elements in $Z_p$. Then $S_A \cap \{0\} = Z_p$. On the
other hand, there is a sequence of $p-2$ non-zero elements of
$Z_p$ such that $S_A \cap \{0\} \neq Z_p$.
\end{fact}
Let $G$ be an abelian group of size $n$ and $q$ be a prime divisor
of $n$. Let $H$ be a subgroup of size $n/q$. A direct corollary of
 Fact \ref{fact1} is the following
 \begin{fact} \label{fact1'}  If $A$ is an incomplete subset of $G$ and
  $S_{A \cap H}= H$, then $H \backslash A$
has at most $q-2$ elements. Consequently $A$ has at most $n/q
+q-2$ elements. \end{fact}

\begin{definition} \label{definition:nice} Let $G$ be an abelian group of size
$n$. A subset $A$ of $G$ is {\it nice} if there is a subgroup $H$
of $G$ such that $|G/H|$ is a prime and $S_{A \cap H} = H$.
\end{definition}
Given a subgroup $H$ in $G$ and an element $a \in G$, we use $a/H$
to represent the coset of $H$ which contains $a$. $a/H$ can be
viewed as an element of the quotient group $G/H$. If $B$ is a
subset of $G$, then $B/H :=\{b/H| b \in B\}$ is a sequence in
$G/H$.
\begin{fact} \label{fact2} If $A$ is a nice incomplete set in a
finite abelian group $G$ of size $n$, then $|A| \le \frac{n}{p} +
p-2 $, where $p$ is the smallest prime divisor of $n$.
Furthermore,  $(A \backslash H) /H$ is an incomplete sequence in
$Z_q = G/H$.
\end{fact}
\begin{proof} (Proof of Fact \ref{fact2}) If $A$ is a nice incomplete set
then $|A| \le \frac{n}{q} + q-2 $, for some prime divisor $q$ of
$n$. On the other hand, it is easy to see that $\frac{n}{q} + q
\le \frac{n}{p} +p$, where $p$ is the smallest prime divisor of
$n$. \end{proof}

Our leading idea  is that  relatively large incomplete sets are
nice. A special case has been verified by   Gao, Hamidoune,
Llad\'o and Serra \cite{GHLS}. Their theorem can be reformulated
in the current setting as follows
\begin{theorem} \label{theorem:GHLS} Let $G$ be an abelian group
of order $n=ph$, where $p \ge 5$ is the smallest prime divisor of
$n$, $h \ge 15p$ is composite. Let $A$ be an incomplete subset of
at least $\frac{n}{p+2} +p$ elements. Then $A$ is nice.
Furthermore, there is a subgroup $H$ of size $n/p$ such that $S_{A
\cap H} = H$.
\end{theorem}

For any positive $\ep \le 1$, define
\begin{equation} \label{Cep} C(\ep):= \sqrt{  \frac{40/\ep^2}{ \log
(2/\ep)}  }
\end{equation}
\noindent and let $n(\ep)$ be the smallest integer $m$ such that
for any $n \ge m$
\begin{equation} \label{Nep} n \ge C(\ep) \sqrt {n \log n} >
\frac {4}{\ep^2}.
\end{equation}

\begin{remark} $n(\ep)$ is relatively small.
One can take, say, $n(\ep) = 500/\ep^4$.
\end{remark}

Now we are ready to state our first theorem.
\begin{theorem} \label{theorem:dir1} Let $\delta$ be a positive
constant at most $1/6$ and  $p_1 \le p_2 \dots \le p_t$, $t \ge
2$, be primes satisfying three conditions
\begin{itemize}
\item $p_2 \ge 3$; \item  $n:=\prod_{i=1}^t p_i \ge n(\delta) $;
\item  $p_1 \le \frac{1}{3C(\delta)} \sqrt {n/ \log n}$.
\end{itemize}
 Let $G$ be an abelian group of order $n$ and  $A$ be an incomplete
subset of $G$ of size at least $(5/6 +\delta)\frac{n}{p_1}$. Then
$A$ is nice and  there is a subgroup $H$ such that $n/|H|$ is one
of the $p_i$,  $|A \backslash H| < 3p_1$ and  $S_{A \cap H} = H$.
\end{theorem}

\begin{remark} Let us have a few comments on this theorem.

\begin{itemize}
\item Using Theorem \ref{theorem:dir1} and Facts  \ref{fact1'}, we
can recover a large portion of Theorem \ref{theorem:old}. To see
this, consider an incomplete set $A$ which does not contain zero.
If $|A| \le n/p_1$, there is nothing to prove. If $|A| \ge n/p_1$,
and $n=|G|$ satisfies the assumptions of Theorem
\ref{theorem:dir1}, apply this theorem to obtain the subgroup $H$.
As $A$ does not contain zero,
 then $|A \cap H| \le |H|-1$. By Facts
\ref{fact1'}, $|A| \le n/q +q-3$, where $q= n/|H|$. But $n/q +q
\le n/p_1 +p_1$, so $|A| \le n/p_1+p_1 -3$.

\item The third assumption $p_1 \le \frac{1}{C} \sqrt {n/ \log n}$
in Theorem \ref{theorem:dir1} can be voided if we assume $t \ge 3$
(i.e., $n/p_1$ is composite) and $n$ sufficiently large. In that
case $p_1 \le n^{1/3} \ll \sqrt{n /\log n}$. It follows  that the
assumptions of Theorem \ref{theorem:dir1} are satisfied whenever
$p_1 \ge 3$, $h$ is composite and $n$ is sufficiently large. Thus,
we have a new proof of Didderich conjecture for sufficiently large
$n$. \item Unlike Theorem \ref{theorem:GHLS}, one cannot conclude
that $H$ has size $n/p_1$. It is  easy to give examples where
$|G|/|H|$ can be any of the $p_i$.

\end{itemize}
\end{remark}

The next question is to find the best lower bound on $|A|$ that
guarantees niceness. Our second theorem gives an almost complete
answer for this question.

\begin{theorem} \label{theorem:dir2} For any positive constant  $\delta$
there is  a positive constant $D(\delta)$ such that the following
holds. Let $p_1 \le p_2 \dots \le p_t$, $t \ge 3$, be primes such
that  $p_1 p_2\le \frac{1}{D(\delta)} \sqrt {n/ \log n}$, where
$n:= \prod_{i=1}^t p_i$. Let $G$ be an abelian group of order $n$
and $A$ be an incomplete subset of $G$ with cardinality at least
$(1+ \delta)\frac{n}{p_1 p_2}$. Then $A$ is nice. Furthermore, the
lower bound $(1+ \delta)\frac{n}{p_1 p_2}$ cannot be replaced by
$\frac{n}{p_1 p_2} + n^{1/4-\alpha}$, for any constant $\alpha$.
\end{theorem}
Finally, let us discuss the case when $G=Z_n$, where $n$ is a
prime. This case has not been covered by the results presented so
far. Olson \cite{O2}, improving upon a result of Erd\H os and
Heilbronn \cite{EH}, shows that $c(Z_n) \le \sqrt{4n-3} +1$. His
bound was improved by  da Silva and Hamidoune \cite{daH}  to
$\sqrt {4n-7}$. As far as characterization results are concerned,
we know of the following two results.
\begin{theorem} \label{theorem:DF}

Let $n$ be a prime and $A$ be an incomplete subset of $\BZ_n$  of
size at least $(2n)^{1/2}$. Then there is some non-zero element $b
\in \BZ_p$ such that
$$\sum_{a \in b  A} \| a\| \le n +O(n^{3/4} \log n). $$
\end{theorem}

\begin{theorem} \label{theorem:NSV} Let $n$ be a prime and $A$ be an incomplete
subset of $\BZ_n$  of size at least $1.99 n^{1/2}$. Then there is
some non-zero element $b \in \BZ_p$ such that
$$\sum_{a \in b  A} \| a\| \le n +O(n^{1/2}). $$
\end{theorem}

Theorem \ref{theorem:DF} is due to  Deshouillers and Freiman
\cite{DF}. Theorem \ref{theorem:NSV} is due to Nguyen, Szemer\'edi
and Vu \cite{NSV}. The error term in this $O(n^{1/2})$ is best
possible, as shown by a  construction in \cite{Des2}.

The rest of the paper is organized as follows. Section 3 contains
the main lemma to the proofs, which states that if $A$ is
sufficiently large, then $S_A$ contains a subgroup of size
comparable to $|A|$. The proofs of the theorems come in Sections 4
and 5. Section 6 is devoted to concluding remarks.

\section {The existence of a large subgroup in $S_A$}
Our key tool is the following statement, which asserts that if $A$
is a sufficiently large subset of $G$, then $S_A$ contains a large
subgroup of $G$. Recall the definition of $C(\ep)$ and $n(\ep)$
from \eqref{Cep} and \eqref{Nep}.

\begin{theorem} \label{theorem:subgroup} Let  $0 < \ep < 1$ be a
constant and  $G$ be an abelian  group of size $n$, where $n \ge
\max \{ \frac{4}{\ep^2}, C(\ep)\sqrt { n \log n } \}$. Let $A$ be
a subset of $G$ with at least $  \max \{ \frac{4}{\ep^2},
C(\ep)\sqrt { n \log n }  \}$ elements. Then $S_A$ contains a
subgroup of size at least $(1-\ep)|A|$.
\end{theorem}
\begin{remark}  The bound $(1-\ep)|A|$ is asymptotically sharp,
 as $A$ itself can be a subgroup. The lower bound $C \sqrt {n
\log n}$ for $|A|$ is sharp up to the logarithmic term. To see
this, consider $G= Z_{p^2} $ and $A=\{0,1 \dots, p \}$. It is
clear that $|A| > p = \sqrt {|G|}$. On the other hand, $S_A$ does
not contain any proper subgroup of $G$. It is interesting to see
whether the $\log$ term can be removed.
\end{remark}
\begin{remark}
The theorem also holds for non-abelian group, see Theorem
\ref{theorem:subgroup1} at the end of this section.
\end{remark}
By definition of $n(\ep)$, if  $n \ge n(\ep)$ then
$$n > C(\ep) \sqrt { n \log n }  >
\frac{4}{\ep^2} . $$ \noindent In this case we have the following
corollary, which is easier to use.
\begin{corollary} \label{cor:subgroup}
 Let  $0 < \ep < 1$ be a
constant and  $G$ be an ablian group of size $n \ge n(\ep)$. Let
$A$ be a subset of $G$ with at least $C(\ep) \sqrt { n \log n } $
elements. Then $S_A$ contains a subgroup of size at least
$(1-\ep)|A|$.
\end{corollary}

To prove Theorem \ref{theorem:subgroup}, we use the following
result of Olson \cite{O1} (see also \cite[Chapter 12]{TVbook}).
Let $A$ be a set and $l$ be a positive integer, we define
$$lA:=\{a_1 + \dots a_l | a_i \in A \}. $$
\noindent Also recall that $<A>$ denotes the subgroup generated by
$A$.
\begin{theorem} \label{theorem:olson}  Let $G$ be finite abelian group, $l$
be a positive integer and $0 \in A$ be a finite subset of $G$.
Then either $lA=<A>$ or $|lA| \ge |A| + (l-1) (\frac{|A|}{2} +1)$.
\end{theorem}

Since $|A| + (l-1) (\frac{|A|}{2} +1) \ge (l+1)|A|/2$, the
following corollary is immediate.
\begin{corollary} Let $G$ be a finite abelian  group, $l$
be a positive integer and $0 \in A$ be a finite subset of $G$ such
that $(l+1)|A| \ge 2|G|$, then $$lA = <A>. $$
\end{corollary}
We also needs the following result of Olson \cite{O1}, which
refines an earlier result of  Szemer\'edi \cite{Sz} (Szemer\'edi
proved the theorem for  an unspecified constant instead of 3).
\begin{theorem} \label{theorem:szemeredi}  Let $G$ be an abelian
 group of order $n$ and $A$ be subset of at least $3 \sqrt n$ elements. Then
 $0\in S_A$.
\end{theorem}

We also need the following simple lemma:
\begin{lemma} \label{lemma:cor:1} Let $G$ be an abelian group and $A$ be a subset of
$G$. Let $l$ be a positive integer and $S$ a subset of $G$ such
that every element of $S$ can be represented as the sum of two
different elements of $A$ in at least $2l-1$ ways (not counting
permutations). Then $l S \subset S_A$.
\end{lemma}

\begin{proof} (Proof of Lemma \ref{lemma:cor:1})
Let $x_1, \dots x_l$ be (not necessarily different) elements of
$S$. We represent $x_1 +\dots + x_l$ as a sum of different
elements of $A$ using the greedy algorithm. To start, represent
$x_1 =a_1 + a_1'$ where $a_1 \neq a_1'$ are different elements of
$A$. Assume that we have represented $x_1= a_1 +a_1', \dots, x_i =
a_i +a_i'$, where $1\le i < l$ and $a_1,a_1', \dots, a_i, a_i'$
are all different. Now look at $x_{i+1}$. Each of the $2i$
elements $a_1,a_1', \dots, a_i, a_i'$ appear in at most one
representation of $x_{i+1}$. Since $x_{i+1}$ has at least $2l-1 >
2i$ representations, we can find a representation $x_{i+1} =
a_{i+1}+  a_{i+1}'$ where both $a_{i+1}$ and $a_{i+1}'$ are
different from $a_1,a_1', \dots, a_i, a_i'$. This concludes the
prof.
\end{proof}

\begin{proof} (Proof of Theorem  \ref{theorem:subgroup})
For each element $x \in G$, let
 $m_x$ be the number of ways to represent $x$ as the sum of two
 different elements of $A$ (not counting permutations). A double
 counting argument gives
 \begin{equation} \label{equ:cor:1-1} \sum_{x \in G} m_x = {|A| \choose 2}. \end{equation}
Notice that $m_x$ is at most $M:=|A|/2$.  Set $K:= 2/\ep$.
 Let $S_j$ be the collection of those
$x$ where $ K^{-j} M < m_x \le K^{-j+1} M$ for $j=1, \dots, j_0$,
where $j_0$ is the largest integer such that $K^{-j_0} M \ge 1$.
Let $S_{j_0+1}$ be the collection of those $x$ where $1 \le m_x
\le K^{-j_0} M$. By the definition of $S_j$
\begin{equation} \label{equ:cor:1-2} \sum_{j=1}^{j_0+1} K^{-j+1} M |S_j |
 \ge  \sum_{x \in G} m_x,
\end{equation}
\noindent which, together with   \eqref{equ:cor:1-1} imply
\begin{equation} \label{equ:cor:1-3} \sum_{j=1}^{j_0+1} K^{-j+1} M |S_j |
  \ge  {|A| \choose 2 }
\end{equation}

Call a set $S_j$ {\it small} ($j=1, \dots, j_0+1$) if it has at
most $(1-\ep)|A|$ elements and {\it large} otherwise. The
contribution from the small $S_j$ on the left hand side is at most
$$\sum_{j=1}^{j_0+1} K^{-j+1} M (1-\ep)|A|  \le \frac{K}{K-1} M(1-\ep)|A|
= (1- \frac{\ep}{2-\ep}) \frac{|A|^2}{2}  $$ \noindent taken into
account the facts that $K= 2/\ep$ and $M=|A|/2$. Since $|A| \ge
\frac{4}{\ep^2}$, we have
$$ (1- \frac{\ep}{2-\ep}) \frac{|A|^2}{2} \le (1 -\ep/2)
\frac{|A|^2} {2} - \frac{|A|}{2}. $$ \noindent From this and
\eqref{equ:cor:1-3}, we have

\begin{equation} \label{equ:cor:1-4}
\sum_{S_j \,\, \hbox{large}}  K^{-j+1} M |S_j | \ge {|A| \choose
2} - (1-\ep/2) \frac{|A|^2}{2} + \frac{|A|}{2} = \ep
\frac{|A|^2}{4}.
\end{equation}
\noindent The bound $|A| \ge C(\ep) \sqrt {n \log n} $ guarantees
that
\begin{equation} \label{equ:cor:1-5} \ep \frac{|A|^2}{4}  \ge 5 K n \log_{2/\ep} n.
\end{equation}
(In fact, $C(\ep)$ is defined so that this inequality holds.)  Set
$l_j:= K^{-j} M$. Since the number of large $S_j$ is at most
$j_0+1 \le \lfloor \log_{2/\ep} |A|/2\rfloor +1$,
\eqref{equ:cor:1-4}, \eqref{equ:cor:1-5} and the pigeon hole
principle imply that there is a large $S_j$ such that
$$ l_j  |S_j| \ge 4 n . $$

Notice that $|S_j| \le |G|= n$. It follows that $ (\lfloor l_j/2
\rfloor +1) |S_j| \ge 2n$. Apply Corollary \ref{cor:subgroup} to
$l := \lfloor l_j/2 \rfloor$ and $S:= S_j \cup \{0\}$, we can
conclude that $lS =<S>$.  On the other hand, by the definition of
$S$
 $$lS= \cup_{i=1}^l i S_j \cup \{0\}. $$
\noindent By Lemma \ref{lemma:cor:1},
 $$i S_j \subset S_A, $$
 \noindent for all $1\le i \le l$. Finally, $0 \in S_A$ by Olson's theorem.
  Thus $S_A$ contains $<A>$, which has at least $(1-\ep)|A|$ elements since $S_j$
is large and $|<S> | \ge |S| \ge |S_j|$. This concludes the proof.
\end{proof}

All the tools used in the proof (Theorems \ref{theorem:olson} and
\ref{theorem:szemeredi}, Lemma \ref{lemma:cor:1}) hold for
non-abelian groups. Thus, Theorem \ref{theorem:subgroup} also
holds for this case. The proof requires only two simple
modifications. First, in Lemma \ref{lemma:cor:1}, $2l-1$ is
replaced by $4l-3$. The reason is that in the proof, each of the
elements $a_1, a_1', \dots, a_i, a_i'$ can now appear in at most 2
representations of $x_{i+1}$. The second is that in the proof of
Theorem \ref{theorem:subgroup}, we need to fix an ordering on the
elements of $G$ and when we consider a sum $x+y$, we always assume
that $x$ precedes $y$ in this ordering. The rest of the proof
remains the same.
\begin{theorem} \label{theorem:subgroup1} For any constant  $0 < \ep < 1$
there are constant $n_1(\ep)$ and $C_1(\ep)$ such that the
following holds. Let  $G$ be a group of size $n$, where $n \ge
n_1(\ep)$. Let $A$ be a subset of $G$ with at least $C_1(\ep)\sqrt
{ n \log n } \}$ elements. Then $S_A$ contains a subgroup of size
at least $(1-\ep)|A|$.
\end{theorem}
The values of $n_1(\ep)$ and $C_1 (\ep)$ might be slightly
different from that of $n(\ep)$ and $C(\ep)$, due to the
modifications.

\section{Proof of Theorem \ref{theorem:dir1}}

\begin{lemma} \label{lemma:dir2} Let $G$ be a finite additive
group and $A$ be a subset of $G$ with cardinality at least
$\lfloor |G|/2 \rfloor +2$. Then $S_A= G$.
\end{lemma}
\begin{proof} (Proof of Lemma \ref{lemma:dir2}) Let $x$ be an
arbitrary element of $G$. There are exactly $\lfloor |G|/2
\rfloor$ (unordered) pairs $(a,b)$ of different elements of $G$
such that $a+b=x$. The claim follows  by the pigeon hole
principle. One can improve the bound slightly but from our point
of view it is not important.
\end{proof}

\begin{proof} (Proof of Theorem \ref{theorem:dir1})
Let $A_1$ be an arbitrary subset of $A$ with cardinality $(1+ 2
\delta)\frac{n}{3p_1}$. By the upper bound on $p_1$,  we can
assume that $|A_1| \ge C(\delta) \sqrt {n \log n}$, which enables
us to apply Corollary \ref{theorem:subgroup} to $A_1$ and obtain a
subgroup $H \subset S_{A_1}$ where
$$|H| \ge (1-\delta)|A_1| = (1-\delta) (1+ 2\delta) \frac{n}{3p_1} >
\frac{n}{3p_1}. $$ The assumption  $p_2 \ge 3$ shows that $H >
\frac{n}{p_1p_2}$. It follows that  $|H| = n/q$ where $q$ is one
of the $p_i$ ( $1\le i \le t$). Furthermore,
$$q < 3 p_1 .$$

  Consider the sequence  $B:= \{a/H | a \in
A \backslash A_1 \}$ in the quotient group $ G/H=Z_q$. If $B$ has
at least $q-1$ non-zero elements, then by Fact \ref{fact1} $S_B$
contains $Z_q \backslash \{0\}$, which implies that
$$G \subset S_{A_1} + S_{A\backslash A_1} \subset S_A $$
\noindent  a contradiction as $A$ is incomplete. Thus, $B$ has at
most $q-2$ non-zero elements. So we can conclude that all but at
most $q-2$ elements of $A \backslash A_1$ lie in $H$. Let $A_2$
denote the set of these elements. We have
$$|A_2| > |A \backslash A_1| -(q-2) \ge |A|-|A_1| - 3p_1 +2 \ge
\Big(\frac{5}{6} +\delta - (1+ 2\delta) \frac{1}{3} \Big)
\frac{n}{p_1} - 3p_1  +2. $$ The right most formula is
$$ (\frac{1}{2} + \frac{1}{3} \delta)) \frac{n}{p_1} -
3p_1 +2 \ge \frac{n}{2p_1} +2.
$$

\noindent since $\frac{1}{3} \delta \frac{n}{p_1} \ge 3p_1$ by the
assumption $p_1 \le \frac{1}{3C(\delta)}  \sqrt {n / \log n}$ and
the definition of $C(\delta)$.

 On the other hand, $|H|$ is at most $ \frac{n}{p_1}$. Thus,
$|A_2| \ge |H|/2 +2$ and so by Lemma \ref{lemma:dir2},
$S_{A_2}=H$. Notice that   $A_2 \subset H \cap A$. Thus $S_{A \cap
H} = H$ which means that  $A$ is nice, completing the proof.
\end{proof}

%By Fact \ref{fact2}, $|A \backslash H| \le q-2 < 3p_1$, so
%$$|H| \ge |A \cap H| \ge (\frac{5}{6} + \delta) \frac{n}{p_1} - 3 p_1 >
%\frac{2}{3} \frac{n}{p_1} $$
%\noindent which implies that $|H| = n/p_1$, proving the theorem.

\section {Proof of Theorem \ref{theorem:dir2} }
 Without loss of generality, we can
 assume that $\delta \le 1/2$ and  $A$ has exactly $(1+\delta) \frac{n}{p_1p_2}$ elements.
 Let $A_1$ be a
subset of $A$ of size $(1+\delta/2) \frac{n}{p_1p_2}$. Setting
$D(\delta)$ sufficiently large, one can assume that $n$ is
sufficiently large and $|A_1| \ge C(\delta/4) \sqrt {n \log n}$
(where $C$ is defined as in \eqref{Cep}), thanks to the assumption
$$p_1p_2 \le \frac{1}{D(\delta)} \sqrt{n /\log n}. $$
\noindent  This enables us to apply Corollary \ref{cor:subgroup}
to $A_1$ and conclude that $S_{A_1}$ contains a subgroup $H$ of
size at least
$$(1-\delta/4) |A_1| = (1-\delta/4) (1+\delta/2)\frac{n}{p_1p_2} >
\frac{n}{p_1p_2}.
$$
The critical point here is that $|H|$ is larger than
$\frac{n}{p_1p_2}$. This forces $|H| = n/q$ where $q$ is one of
the primes $p_i$. It would be  easy to finish the proof now if
$A$ had at least $(2+\delta) \frac{n}{p_1p_2}$ (instead of only
$(1+\delta) \frac{n}{p_1p_2}$) elements. The reason is that in
this case we still have $(1+\delta/2)\frac{n}{p_1p_2}$ elements
outside $A_1$ to play with. Arguing as in the previous proof, we
can show that most of these elements should be in $H$ and span it.
As we lack these extra elements, we need an additional trick that
helps us to show that actually most elements of $A_1$ are already
in $H$. The heart of this trick is Lemma \ref{lemma:random} below.
Before presenting the lemma, let us make some observations. Set
$A_2 := A \backslash A_1$. As $A_1$ was chosen arbitrarily, $A_2$
is an arbitrary subset of $A$ with  $\frac{\delta n}{2 p_1p_1}$
elements. Since $A$ is incomplete, $|A_2 \backslash H| \le q-2$,
where $q=|G|/|H| \le p_1p_2$. By setting $D(\delta)$ sufficiently
large, we can assume
$$p_1p_2 \le \frac{\delta^2}{20}  \frac{n}{p_1p_2} =
\frac{\delta}{10} |A_2| $$

\noindent which implies
$$|A_2 \cap H| \ge  |A_2| - p_1p_2 \ge (1-\delta/10) |A_2|. $$

To summarize, $A$ has the property that for any subset $A_2$ of
size $\frac{\delta n}{2 p_1p_2} = \frac{\delta}{2(1+\delta)} |A|$,
there is a maximal subgroup $H$ of $G$ such that $|A_2 \cap H| \ge
(1-\delta/10) |A_2|$.

\begin{lemma} \label{lemma:random}
Let $S$ be a subset of $G$ of size $(1+\delta) \frac{n}{p_1p_2}$
such that no maximal subgroup of $G$ contains $(1-\delta/2)$
fraction of $S$. Then there is a subset $S' \subset S$ of size
$\frac{\delta}{2(1+\delta)} |S|$ such that no maximal subgroup of
$G$ contains $(1-\delta/10  )$ fraction of $S'$.
\end{lemma}

Assuming the lemma for a moment, we can conclude the proof as
follows.  By the lemma and its preceding paragraph, there  is a
maximal subgroup $H$ such that
$$|H \cap A| \ge (1- \delta/2) |A| = (1-\delta/2)(1+\delta) \frac{n}{p_1p_2} \ge
(1+\delta/4) \frac{n}{p_1p_2} $$ as $\delta \le 1/2$. Since $|H|
\le n/p_1$ and the smallest prime divisor $p'$ of $H$ is either
$p_1$ or $p_2$, it is easy to verify that
$$|H \cap A| \ge \frac{|H|} {p'} + p' . $$
\noindent Thus we can apply Theorem \ref{theorem:old} or Theorem
\ref{theorem:dir1}  for $H$ and $A \cap H$ to deduce  that $A \cap
H$ is complete in $H$. Therefore,  $S_{A \cap H} = H$ and $A$ is
nice.

Now we prove  Lemma \ref{lemma:random}, using a probabilistic
argument.
\begin{proof} (Proof of Lemma \ref{lemma:random}) Set $s:=
|S|= (1+\delta) \frac{n}{p_1p_2}$ and $\ep:= \delta/10$.  Consider
a random subset $S_1$ of $S$ obtained by  selecting each element
$a \in S$ to be in $S_1$ with probability $\rho:= (1+2\ep)
\frac{\delta}{2(1+\delta)^2}$, independently. Let $H$ be  a
subgroup of $G$. By linearity of expectation and the assumption of
the lemma , we have
$$\E (|H \cap S_1|) = \rho |H \cap S| \le \rho (1-\delta/2) s = \rho (1-5\ep) s. $$
\noindent On the other hand, $\E(|S_1| ) = \rho s$. Both $H \cap
S_1$ and $S_1$ have binomial distribution.  By property of the
binomial distribution, there is a positive constant $c_0$
depending only on $\ep$  such that  with probability at least $1-
\exp(-c_0\rho s)$
\begin{equation} \label{S1-1}  (1-\ep) \rho s \le |S_1| \le (1+\ep) \rho s
. \end{equation} \noindent and
\begin{equation} \label{S1-2}  |H \cap S_1| \le (1+\ep)  \rho (1-5\ep) s.
\end{equation}

\noindent  It is well known (and easy to prove) that  the number
of maximal subgroups  of $G$ is at most $|G|= n$. If $D(\delta)$
(and so $n$) is sufficiently large, then
 $$2n \le \exp( c_0\rho s ). $$
 \noindent Thus, we can use the union bound to conclude that there exists a set
  $S_1$ such that \eqref{S1-1} holds and \eqref{S1-2} holds
simultaneously for every maximal subgroup $H$. Let $S'$ be any
subset of $S_1$ of size $\frac{\delta n}{2 n_1n_2}=
\frac{1}{1+2\ep} \rho s$. For any maximal subgroup $H$
$$|S' \cap H|/ |S'|  \le |S_1 \cap H| /|S'| \le \frac{(1+\ep)(1-5\ep)\rho
s}{(\frac{1}{1+2\ep} \rho s} < (1-\ep) =(1-\delta/10). $$
\noindent This concludes the proof of the lemma. \end{proof}

The following example shows  that the lower bound $(1+\delta)
\frac{n}{p_1p_2}$ cannot be reduced to $\frac{n}{p_1p_2} +
n^{1/4-\alpha}$, for any fixed $\alpha$.

{\it Example.} Take $n:= p^2 q$ where $1 < p < q $ are large
primes. Consider $G= Z_{p^2} \oplus Z_q$. Given any $\delta >0$
and any function $D(\delta)$,  by choosing $q$ properly $p$ we can
guarantee that
$$ n^{1/2-\alpha} \le p^2 \le \frac{1}{D(\delta)} \sqrt {n/ \log n}. $$

 We write an element $a \in G$ as $a= (x,y)$
where $x \in Z_{p^2} $ and $y \in Z_q$. Let $m$ be the largest
integer such that $\sum_{i=0}^m i < p^2-1$. Set
$$A:= \{ (x, 0) | 0 \le x \le m \} \cup \{ (0, y)|
0\le y \le q -1 \} .$$ It is easy to show that $A$ is  incomplete
and not nice, thanks to the fact that $\sum_{i=0}^m i < p^2 -1$.
On the other hand,
$$|A| = m + q = m + \frac{n}{p^2} = m + \frac{n}{p_1p_2} \ge n^{1/4 -\alpha}
 + \frac{n}{p_1p_2} . $$
The proof of the theorem is complete.

\section{Concluding remarks}

One can use the additional trick in the proof of Theorem
\ref{theorem:dir2} to improve upon the constant  $(5/6+\delta)$ in
Theorem \ref{theorem:dir1}. However, this requires some
modification on the assumptions. We prefer to present Theorem
\ref{theorem:dir1} in the simplest way in order to illustrate the
ideas.

One can also use the method presented here to study incomplete
sets with size less than $\frac{n}{p_1p_2}$. However, the
characterization obtained in this case is more technical and less
appealing.

\end{document}